\documentclass{amsproc}

\newtheorem{theorem}{Theorem}[section]

\newtheorem{proposition}[theorem]{proposition}

\theoremstyle{definition}
\newtheorem{definition}[theorem]{Definition}

\theoremstyle{remark}
\newtheorem{remark}[theorem]{Remark}

\numberwithin{equation}{section}

\usepackage{color}
\newcommand{\vphi}{\vcg{\varphi}}
\newcommand{\kappad}{\kappa_{\delta}}

\newcommand{\N}{\mathbf{N}}
\newcommand{\dx}{{\, \rm d}x}

\newcommand{\Div}{{\rm div}\,}
\newcommand{\Ov}[1]{\overline{#1}}

\newcommand{\Un}[1]{\underline{#1}}
\newcommand{\bo}{| _{\partial\Omega}}

\newcommand{\vy}{Y_{k}}
\newcommand{\vf}{{\vc F}}

\newcommand{\vw}{\omega}

\newcommand{\vr}{\varrho}

\newcommand{\vp}{\varphi}

\newcommand{\vrd}{\vr_\delta}
\newcommand{\vcYd}{\vec{Y}_\delta}

\newcommand{\vtd}{\vt_\delta}
\newcommand{\vud}{\vu_\delta}
\newcommand{\vt}{\vartheta}

\newcommand{\vu}{\vc{u}}

\newcommand{\vc}[1]{{\bf #1}}
\newcommand{\vcg}[1]{{\pmb #1}}
\newcommand{\F}[1]{$\mathbb{#1}$}
\newcommand{\Grad}{\nabla}

\newcommand{\tn}[1]{\mbox {\F #1}}

\newcommand{\dS}{{\rm d}S}

\newcommand{\lr}[1]{\left( #1 \right)}
\newcommand{\intO}[1]{\int_{\Omega} #1\dx}
\newcommand{\intpO}[1]{\int_{\partial\Omega} #1 \, {\rm d} {S}}

\newcommand{\intOB}[1]{\int_{\Omega} \left( #1 \right) \ \dx}

\newcommand{\D}{{\vc{D}}}
\newcommand{\ep}{\varepsilon}

\newcommand{\R}{\mathbb{R}}
\newcommand{\sumkN}[1]{\sum_{k=1}^n #1}
\newcommand{\sumlN}[1]{\sum_{l=1}^n #1}

\begin{document}

\title[Chemically reacting heat conducting compressible mixture]{On steady solutions to a model of chemically reacting heat conducting compressible mixture with slip boundary conditions}

\author{Tomasz Piasecki}
\address{Institute of Applied Mathematics and Mechanics, University of Warsaw, ul. Banacha 2, 02-097 Warsaw, Poland}
\email{tpiasecki@mimuw.edu.pl}
\thanks{The work of the first author (TP) was supported  by the Polish NCN grant No. UMO-2014/14/M/ST1/00108}

\author{Milan Pokorn\'y}
\address{Charles University, Faculty of Mathematics and Physics, Sokolovsk\'a 83, 18675 Praha 8, Czech Republic}
\email{pokorny@karlin.mff.cuni.cz}
\thanks{The work of the second author (MP) was supported by the Czech Science Foundation (grant No. 16-03230S)}

\date{June 30, 2017 and, in revised form, }



\begin{abstract}
We consider a model of chemically reacting heat conducting compressible mixture. 
We investigate the corresponding system of partial differential equations in the steady regime 
with slip boundary conditions for the velocity and, in dependence on the model parameters, 
we establish existence of either weak or variational entropy solutions. 
The results extend the range of parameters for which the existence of weak solutions is known in the case of homogeneous Dirichlet boundary conditions for the velocity.
\end{abstract}

\maketitle

\section{Introduction}

The aim of this note is to extend results on a model of chemically reacting heat conducting compressible gaseous mixture based on the model considered e.g. in \cite{Gi}. Results dealing with steady solutions for this model appeared recently in \cite{GPZ} and \cite{PiPo} for the case of Dirichlet boundary conditions; see also \cite{Za} which can be considered actually as a first result in this direction (dealing, however, with a slightly simplified model).     

The common feature of these three papers is the fact that the 
weak solutions (and also variational entropy solutions in \cite{PiPo}) were obtained for any relatively rough data, without any assumption on their size or on the distance to a known (possibly regular) solution. 

This paper is devoted to the proof of existence of weak and variational entropy solutions to the model introduced in \cite{GPZ}. Due to the slip boundary conditions for the velocity we are able to extend the range of parameters for which the weak solutions exist.  This corresponds to the fact which has been observed several times for the compressible Navier--Stokes or the Navier--Stokes--Fourier system; the slip boundary conditions allow for better density (and sometimes also velocity) estimates which leads to stronger results in this case, see e.g. \cite{MuPo1}, \cite{PoMu}, \cite{JN}, \cite{MuPo2}, \cite{JNP} or also \cite{MPZ_Handbook}.

In what follows, we use standard notation for the Lebesgue, Sobolev and other standard function spaces as well as for norms in these spaces. The scalar valued functions will be denoted by the standard font (e.g., $\vr$ and $\vt$ for the density and temperature, respectively), the vector valued functions will be printed in bold face (e.g., $\vu$ for the velocity) and the tensor valued functions using a special font (e.g., $\tn{S}$ for the viscous part of the stress tensor). The generic constants are denoted by $C$ and their value may change from line to line or even in the same formula.

\subsection{The model}

We consider the following system of partial differential equations
\begin{equation}\label{1}
\begin{array}{c}
\Div (\vr \vu) = 0,\\
\Div (\vr \vu \otimes \vu) - \Div \tn{S} + \Grad \pi =\vr \vc{f},\\
\Div (\vr E\vu )+\Div(\pi\vu) +\Div\bf{Q}- \Div (\tn{S}\vu)=\vr\vc{f}\cdot\vu,\\
\Div (\vr Y_k \vu)+ \Div \vf_{k}  =  m_k\vw_{k},\quad k\in \{1,\ldots,n\},
\end{array}
\end{equation}
where the unknown quantities are the total density $\vr$, the velocity field $\vu$, the temperature $\vt$ (appearing in (\ref{1}) implicitly, see below) and the mass fractions $Y_k = \vr_k/\vr$, where $\{\vr_k\}_{k=1}^n$ are the densities of the constituents. As $\sumkN \vr_k = \vr$, we have $\sumkN Y_k =1$. The other functions, i.e. the stress tensor $\tn{S}$, the pressure $\pi$, the total energy $E$, the heat flux $\bf{Q}$, the diffusion fluxes $\vf_k$ and the molar production rates $\vw_k$ are given functions of these unknows and will be introduced below. Furthermore, $\vc{f}$ is the given field of the external forces (e.g., the gravity force) and $m_k$ denotes the molar masses of the $k$th constituent, $k=1,2,\dots, n$.  
 
System (\ref{1}) is completed by the boundary conditions on $\partial \Omega$ 
\begin{equation} \label{2}
\begin{array}{rcl}
 \vf_{k}\cdot\vc{n}&=&0, \\
-\vc{Q}\cdot\vc{n}+L(\vt-\vt_{0})&=&0,
\end{array}
\end{equation}
and for the velocity we assume  the Navier boundary condition (the slip b.c.)
\begin{equation} \label{4}
\vu \cdot \vc{n} = 0, \qquad (\tn{S} \vc{n} + f \vu)\times \vc{n} = \vc{0}.
\end{equation}
Above, the boundary condition for the temperature means that the heat flux is proportional to the difference of the temperature inside and outside the boundary. The coefficient $f$ (assumed to be constant in what follows) denotes the friction. 
We also prescribe the total mass
\begin{equation} \label{5}
\intO {\vr} = M>0.
\end{equation}

\subsubsection{The stress tensor and the pressure}

We assume the stress tensor $\tn{S}$ to be a given linear function of the symmetric part of the velocity gradient
\begin{equation} \label{6}
{\tn S} = {\tn S}(\vt, \widetilde{\tn{D}}(\vu))= \mu\Big[\Grad \vu+(\Grad \vu)^{T}-\frac{2}{3}\Div \vu \tn{I}\Big]+\nu(\Div \vu)\tn{I},
\end{equation}
where $\widetilde{\tn{D}}(\vu) = \frac 12 (\Grad \vu+(\Grad \vu)^{T})$, the coefficients $\mu=\mu(\vt)>0$ (Lipschitz continuous in $\R^+$) $\nu=\nu(\vt)\geq 0$ (continuous in $\R^+$),  are the shear and bulk viscosity coefficients, respectively. We assume
\begin{equation} \label{7}
\Un{\mu}(1+\vt)\leq\mu(\vt)\leq\Ov{\mu}(1+\vt),\quad 0\leq\nu(\vt)\leq\Ov{\nu}(1+\vt)
\end{equation}
for positive constants $\underline{\mu},\overline{\mu},\overline{\nu}$.
Furthermore, $\tn{I}$ is the identity matrix.

The pressure
\begin{equation}\label{8}
\pi =\pi(\vr,\vt)=\pi_{c}(\vr)+\pi_{m}(\vr,\vt),
\end{equation}
where the cold pressure is assumed in the form
\begin{equation} \label{9}
\pi_{c}=\vr^{\gamma}, \quad \gamma>1.
\end{equation}
A more general pressure form can be assumed, as in the case of the steady compressible Navier--Stokes--Fourier system, see e.g. \cite{NoPo_JDE} or \cite{FPT} in the context of the chemically reacting flows. We, however, prefer to keep its form as simple as possible.
The molecular pressure $\pi_{m}$, according to the Boyle law, satisfies 
\begin{equation}\label{9a}
\pi_{m}=\pi_m(\vr,\vt) = \sumkN p_k(\vr,\vt) = \sumkN\frac{\vr Y_k}{m_k}\vt,
\end{equation}
where, for simplicity, the gas constant is taken to be equal to one.

\subsubsection{The energy and the heat flux}

The specific total energy $E$ is a sum of the specific kinetic and specific internal energies (we denote $\vec Y = (Y_1,Y_2,\dots,Y_n)$)
\begin{equation} \label{10}
E= E(\vr,\vu,\vt,\vec Y)=\frac{1}{2}|\vu|^{2}+e(\vr,\vt,\vec Y).
\end{equation}
Due to the form of the pressure the internal energy consists of two components 
\begin{equation} \label{11}
e=e_{c}(\vr)+e_{m}(\vt,\vec Y),
\end{equation}
where the cold energy $e_c$ and the molecular internal energy $e_m$ are given by
\begin{equation} \label{12}
e_{c}=\frac{1}{\gamma-1}\vr^{\gamma-1},\qquad\qquad e_m= \sumkN Y_ke_k=\vt\sumkN c_{vk}Y_k.
\end{equation}
Above, $c_{vk}$ are the constant-volume specific heats and can be different for different species. They are related to the 
constant-pressure specific heats by
\begin{equation}\label{13}
c_{pk}=c_{vk}+\frac {1}{m_k}
\end{equation} 
and both $c_{vk}$ and $c_{pk}$ are assumed to be constant. 

The heat flux $\bf{Q}$ consists of two terms 
\begin{equation}\label{13a}
\vc{Q}=\sumkN h_k \vf_{k}+\vc{q},
\end{equation}
where the first term represents transfer of energy due to the species molecular diffusion (and $h_k$, defined below, are the enthalpies) and the second one the Fourier law
\begin{equation}\label{13b}
\vc{q}=-\kappa\Grad\vt.
\end{equation}
The coefficient $\kappa=\kappa(\vt)$ is the thermal conductivity coefficient and we assume 
\begin{equation} \label{13c}
\underline{\kappa}(1+\vt^{m})\leq\kappa(\vt)\leq\Ov{\kappa}(1+\vt^{m})
\end{equation}
for some constants $m,\underline{\kappa},\overline{\kappa}>0$.

\subsubsection{Diffusion flux and species production rates}

The form of the diffusion flux is the most important part modeling the interaction between the species. Following \cite{Gi} we assume that
\begin{equation} \label{14}
\vf_k = -\sumlN C_{kl} \vc{d}_l, \quad k=1,2,\dots, n,
\end{equation}
where 
\begin{equation}
\label{15}
\vc{d}_k = \Grad \Big(\frac{p_k}{\pi_m}\Big) + \Big( \frac{p_k}{\pi_m} - \frac{\vr_k}{\vr}\Big) \Grad \log \pi_m
=\frac{\nabla p_k}{\pi_m}-Y_k\frac{\nabla \pi_m}{\pi_m}.
\end{equation}
Furthermore, we introduce another matrix $\tn{D}$ 
$$
C_{kl} = Y_k D_{kl},
$$
where the diffusion matrix $\tn{D}=\tn{D}(\vt,\vec{Y})$ has the following properties
\begin{equation}\label{pr15}
	\begin{gathered}
		\tn{D}=\tn{D}^T,\quad
		N(\tn{D})=\R  \vec{Y},\quad
		R(\tn{D})={\vec{U}}^{\bot},\\
		\tn{D} \quad\text{ is positive semidefinite over } \R^n,
	\end{gathered}
	\end{equation}
with $\vec{Y}=(Y_1,\ldots,Y_n)^T>0$ and $\vec{U} = (1,\dots, 1)^T$. 
Above $N(\tn{D})$ denotes the nullspace of matrix $\tn{D}$, $R(\tn{D})$ denotes its range, and ${\vec{U}}^{\bot}$ denotes the orthogonal complement of $\vec{U}$. Moreover, the matrix $\tn{D}$ is positively definite over $\vec{U}^{\bot}$ and there exists $\delta>0$ such that
\begin{equation} \label{16}
\delta \langle \tn{Y}^{-1}\vec{x},\vec{x}\rangle \leq \langle \tn{D}\vec{x},\vec{x}\rangle \quad \forall \vec{x} \in \vec{U}^{\bot},
\end{equation}
where $\tn{Y}= {\rm diag}\, (Y_1,\dots,Y_n)$ and $\langle \cdot,\cdot\rangle$ denotes the scalar product in $\R^n$. 

Furthermore,  $D_{ij}$ are differentiable functions of $\vt,Y_1,\ldots,Y_n$ for any $i,j\in\{1,\ldots,n\}$ such that 
\begin{equation} \label{1.21a}
|Y_iD_{ij}(\vt,\vec{Y})| \leq C(\vec{Y}) (1+\vt^a)
\end{equation}
for some $a\geq 0$, and $C(\vec{Y})$ is bounded in $[0,1]^n$. Finally,
\begin{equation} \label{17}
\sumkN \vf_k=\vc{0}.
\end{equation} 
Note that we can also consider the Fick law in the form
\begin{equation} \label{1.22a}
\vf_k = D(\vt,\vec Y) \Grad Y_k, \quad k=1,2,\dots, n,
\end{equation}
and the function $D(\cdot,\cdot)$ is a differentiable function fulfilling a similar estimate as (\ref{1.21a}), i.e.
\begin{equation} \label{1.22b}
0<D_0 \leq D(\vt,\vec Y) \leq C(\vec Y) (1+\vt^a)
\end{equation}
with $a\geq 0$ and $C(\cdot)$ bounded in $[0,1]^n$. Condition (\ref{17}) is indeed fulfilled provided $\sumkN Y_k =1$.

Concerning the species production rates, we assume that $\{\omega_k\}_{k=1}^n$ are differentiable functions of $\vr,\vt,\vec{Y}$ which are bounded, and
such that
\begin{equation} \label{18}
\omega_k \geq - CY_k^r \quad \textrm{for some} \quad C,r>0,
\end{equation}
which means that a species cannot decrease faster 
than proportionally to some positive power of its fraction (a possible natural choice is $r=1$). Moreover, this condition clearly implies the compatibility condition $\omega_k \geq 0$ if $Y_k=0$. Furthermore,
\begin{equation} \label{19}
\sumkN m_k \omega_k =0.
\end{equation}
Note that e.g. in \cite{FPT}, instead of $m_k \omega_k$ it is assumed that the species production rate is modeled as $\vr m_k\omega_k$. We may easily treat here this version and in a sense (see comments to the entropy inequality below) it is in fact simpler.

\subsubsection{Entropy and other thermodynamic potentials; entropy production rate}

Since our main thermodynamic quantities are the internal energy and the pressure, the other thermodynamic potentials are assumed to be given functions of them. In what follows, we assume the thermodynamics connected with a mixture of ideal gases with addition of the cold pressure term and the corresponding term in the internal energy. 

The specific enthalpy (of each constituent) has the form
\begin{equation} \label{20}
h_k = c_{pk} \vt, \qquad h = \sumkN Y_k h_k, 
\end{equation}
where $c_{pk}$ fulfills (\ref{13}). The specific entropy
\begin{equation} \label{21}
s_k = c_{vk} \log \vt - \frac{1}{m_k} \log \Big(\frac{\vr Y_k}{m_k}\Big), \qquad s= \sumkN Y_k s_k,
\end{equation}
and the Gibbs function (Gibbs free energy) 
\begin{equation} \label{22}
g_k = h_k -\vt s_k, \qquad g = \sumkN Y_k g_k.
\end{equation}
Moreover, the Gibbs formula has the form
\begin{equation}\label{23}
\vt \D s=\D e+\pi\D\left({\frac {1}{\vr}}\right)-\sumkN g_{k}\D \vy.
\end{equation}
Using (\ref{23}) it is possible to derive an equation for the specific entropy $s$
\begin{equation}\label{24}
\Div(\vr s\vu)+\Div\left( \frac{\vc{Q}}{\vt}-\sumkN \frac{g_{k}}{\vt}\vf_{k}\right)=\sigma,
\end{equation}
where the entropy production rate
\begin{equation} \label{25}
\sigma=\frac{{\tn S}:\Grad\vu}{\vt}-{\frac{\vc{Q}\cdot\Grad\vt}{\vt^{2}}}-\sumkN\vf_{k}\cdot\Grad\left({\frac{g_{k}} {\vt}}\right)-\frac{\sumkN m_k g_{k}\vw_{k}}{\vt}.
\end{equation}
Note that the entropy production rate can be expressed in the form
\begin{equation} \label{26}
\sigma=\frac{{\tn S}:\Grad\vu}{\vt}+{\frac{\kappa |\Grad\vt|^2}{\vt^{2}}}-\sumkN\frac{\vf_{k}}{m_k}\cdot\Grad\log p_k-\frac{\sumkN m_k g_{k}\vw_{k}}{\vt}.
\end{equation}
Then we easily see that the first two terms are non-negative due to the form of the stress tensor and the positivity of $\kappa$. Moreover, we assume that
\begin{equation} \label{26a}
\sumkN m_k g_k \omega_k \leq 0,
\end{equation}
which implies that also the fourth term is non-negative. Finally,
$$
\begin{aligned}
-\sumkN\frac{\vf_{k}}{m_k}\cdot\Grad\log p_k &= \sum_{k,l=1}^n \frac{Y_k D_{kl}}{m_k} \Big[\Grad \Big(\frac{p_l}{\pi_m}\Big) + \Big(\frac{p_l}{\pi_m}-\frac{\vr_l}{\vr}\Big)\Grad \log \pi_m \Big] \Grad \log p_k \\
&= \sum_{k,l=1}^n \frac{Y_k D_{kl}}{m_k} \Big[ \frac{\Grad p_l}{\pi_m} -\frac{\vr_l}{\vr} \frac{\Grad \pi_m}{\pi_m} \Big] \frac{\Grad p_k}{p_k}\\
&= \frac{\pi_m}{\vr\vt}\sum_{k,l=1}^n D_{kl} \Big[ \frac{\Grad p_l}{\pi_m} -Y_l \frac{\Grad \pi_m}{\pi_m} \Big] \Big[ \frac{\Grad p_k}{\pi_m} -Y_k \frac{\Grad \pi_m}{\pi_m} \Big]\\
+  \frac{\pi_m}{\vr\vt}& \sum_{k,l=1}^n D_{kl} \Big[ \frac{\Grad p_l}{\pi_m} -Y_l \frac{\Grad \pi_m}{\pi_m} \Big] Y_k \frac{\Grad \pi_m}{\pi_m} \geq 0
\end{aligned}
$$
due to the properties of the matrix $\tn{D}$, as $\sumkN \vf_k = \vc{0}$ implies
\begin{equation} \label{27}
\sum_{k=1}^n Y_k D_{kl} = 0 \quad \forall l=1,2,\dots, n.
\end{equation}
Notice in particular that Fick law \eqref{1.22a} is not a special case of (\ref{14})--(\ref{pr15}).  

\subsection{Formulation of the problem}

We now formulate the problem we treat in this paper.

\subsubsection{Non-diagonal diffusion matrix}

Unfortunately, the general non-diagonal form of the diffusion flux is too complex to be considered in the full generality. In particular, the above deduced lower bound of the corresponding term in the entropy production rate does not allow us to control the gradient of the mass fractions. A certain attempt has been done in the evolutionary case, see \cite{MPZ}, however, it leads to the necessity to control the (total) density gradient, which can be obtained for the fluids with density dependent viscosities satisfying the Bresch--Desjardins identity. The same idea does not work in the steady problem and therefore we must restrict ourselves (cf. \cite{GPZ} or \cite{PiPo}) to the case when all molar masses are comparable. We therefore assume that $m_1=m_2=\dots=m_n$ and without loss of generality we set this common value to be equal to one.

Then $\pi_m = \sumkN \vr_k \vt = \vr\vt$, $\frac{\Grad p_l}{\pi_m} -Y_l \frac{\Grad \pi_m}{\pi_m} = \Grad Y_l + Y_l \frac{\Grad (\vr\vt)}{\vr\vt} - Y_l \frac{\Grad (\vr\vt)}{\vr\vt}= \Grad Y_l$. 
Therefore
$$
\vf_k=-\sum_{l=1}^n Y_k D_{kl}\nabla Y_l
$$
and
$$
-\sumkN\frac{\vf_{k}}{m_k}\cdot\Grad\log p_k = -\sumkN\vf_{k}\cdot\Big(\frac{\Grad Y_k}{Y_k}+\frac{\nabla(\vr\vt)}{\vr\vt}\Big)= \sumkN D_{kl} \Grad Y_l \Grad Y_k \geq c |\Grad \vec{Y}|^2,
$$
provided $\vec{Y} \geq 0$ and $\sumkN Y_k =1$. Exactly this estimate allows to obtain the existence of a solution in this case. We may therefore consider

\bigskip

\noindent {\bf Problem P} (Non-diagonal diffusion matrix) 

We consider system (\ref{1}) with boundary conditions (\ref{2}), (\ref{4}), given total mass (\ref{5}), and (\ref{6})--(\ref{17}), (\ref{18})--(\ref{25}) with equal molar masses $m_1=m_2=\dots=m_n =1$.

\bigskip

\subsubsection{Fick's law}

A similar estimate as above we get also in the case of the Fick law (\ref{1.22a}). We have for the molar masses being the same (and for notational simplicity, equal to 1)
\begin{multline*}
-\sumkN\frac{\vf_{k}}{m_k}\cdot\Grad\log p_k = -\sumkN\vf_{k}\cdot\Grad\log p_k  \\ =\sumkN \frac{D(\vt,\vec Y)}{Y_k} \nabla Y_k\cdot \nabla Y_k +  \sumkN D(\vt,\vec Y) \nabla Y_k \cdot \frac{\nabla (\vr\vt)}{\vr \vt} \geq D_0 |\nabla \vec Y|^2.
\end{multline*}
Note that in the case of the same molar masses the Fick law behaves exactly in the same way as Problem P for the same molar masses and therefore we do not consider it separately.

\section{Definitions of solutions. Existence results}

Problem P with the Dirichlet boundary conditions for the velocity has been studied in \cite{GPZ} and \cite{PiPo}, we therefore present in analogy to this case the definitions of weak and variational entropy solutions, in the spirit of paper \cite{JNP}. We introduce
$$
C^1_{\vc{n}} (\Omega) = \{\vc{w} \in C^1(\Ov{\Omega}); \vc{w} \cdot \vc{n}=0 \text{ on } \partial \Omega\}.   
$$

\begin{definition}\label{d1}
We say the set of functions $(\vr,\vu,\vt, \vec{Y})$ is a weak solution to system
(\ref{1}) with boundary conditions (\ref{2}), (\ref{4}), given total mass (\ref{5}), and (\ref{6})--(\ref{17}), (\ref{18})--(\ref{25})  with equal molar masses $m_1=m_2=\dots=m_n =1$
provided 
\begin{itemize}
\item
$\vr \geq 0$ a.e. in $\Omega$, $\vr \in L^{6\gamma/5}(\Omega)$, $\int_{\Omega} \vr\dx=M$ 

\item
$\vu \in W^{1,2}(\Omega)$, $\vu\cdot \vc{n} = 0$ a.e. on $\partial \Omega$, $\vr |\vu|$ and $\vr |\vu|^2 \in L^{\frac{6}{5}}(\Omega)$

\item 
$\vt \in W^{1,2}(\Omega) \cap L^{3m}(\Omega)$, $\vr \vt, \vr\vt|\vu|, {\tn S}\vu, \kappa|\nabla \vt| \in L^1(\Omega)$

\item
$\vec{Y}\in W^{1,2}(\Omega)$, $Y_k \geq 0$ a.e. in $\Omega$, $\sumkN Y_k = 1$ a.e. in $\Omega$, $\vf_k\cdot \vc{n}|_{\partial \Omega}=0$
\end{itemize}
and the following integral equalities hold\\
$\bullet$ the weak formulation of the continuity equation
\begin{equation}\label{weak_cont}
\intO{\vr \vu\cdot\Grad\psi} = 0
\end{equation}
holds for any test function $\psi\in C^{1}(\Ov{\Omega})$;\\
$\bullet$ the weak formulation of the momentum equation
\begin{equation} \label{weak_mom}
-\intO{\big(\vr\left(\vu\otimes\vu\right):\Grad\vcg{\vp}-\tn{S}:\Grad \vcg{\vp}\big)}+ f\intpO{\vu \cdot \vcg{\vp}}-\intO {\pi \Div\vcg{\vp}}=\intO{\vr\vc{f}\cdot\vcg{\vp}}
\end{equation}
holds for any test function $\vcg{\vp}\in  C^{1}_{\vc{n}}(\Omega)$;\\
$\bullet$ the weak formulation of the species equations
\begin{equation}\label{weak_spe}
-\intO{ Y_{k}\vr\vu\cdot\Grad\psi}-\intO{\vf_k\cdot\Grad\psi}=\intO{\vw_{k}\psi}
\end{equation}
holds for any test function $\psi\in C^{1}(\Ov{\Omega})$ and for all $k=1,\ldots,n$;\\
$\bullet$ the weak formulation of the total energy balance
\begin{equation} \label{weak_ene}
\begin{aligned}
&-\intO{\left(\frac{1}{2}\vr|\vu|^{2}+\vr e\right)\vu\cdot\Grad\psi}+\intO{\kappa\Grad\vt\cdot\Grad\psi}
-\intO{\left(\sumkN{h_{k}\vf_{k}}\right)\cdot\Grad\psi}\\
&=\intO{\vr\vc{f}\cdot\vu\psi}-\intO{\lr{\tn{S}\vu}\cdot\Grad\psi}+\intO{\pi\vu\cdot\Grad\psi}\\
&-\intpO{L(\vt-\vt_{0})\psi} - f\intpO{|\vu|^2 \psi}
\end{aligned}
\end{equation}
holds for any test function $\psi\in  C^{1}(\Ov{\Omega})$.
\end{definition}

Indeed, the total energy balance which contains the term behaving as $\vr|\vu|^3$ limits the range for $\gamma$ and $m$ for which we are able to prove existence of a weak solution. Following a similar situation for the compressible Navier--Stokes--Fourier system (both steady and evolutionary, see \cite{FeNo} or \cite{MPZ_Handbook}) we introduce another type of solution, where the total energy balance is replaced by the entropy inequality. 

\begin{definition} \label{d2}
We say the set of functions $(\vr,\vu,\vt, \vec{Y})$ is a variational entropy solution to problem 
(\ref{1}) with boundary conditions (\ref{2}), (\ref{4}), given total mass (\ref{5}), and (\ref{6})--(\ref{17}), (\ref{18})--(\ref{25})  with equal molar masses $m_1=m_2=\dots=m_n =1$ 
provided
\begin{itemize}
\item
$\vr \geq 0$ a.e. in $\Omega$, $\vr \in L^{s\gamma}(\Omega)$ for some $s>1$, $\int_{\Omega} \vr\dx=M$ 

\item
$\vu \in W^{1,2}(\Omega)$, $\vu\cdot \vc{n} = 0$ a.e. on $\partial \Omega$, $\vr \vu \in L^{\frac{6}{5}}(\Omega)$

\item
$\vt \in W^{1,r}(\Omega) \cap L^{3m}(\Omega)$, $r>1$, $\vr \vt, \tn {S}:\frac{\nabla \vu}{\vt}, \kappa\frac{|\nabla \vt|^2}{\vt^2}, \kappa\frac{\nabla \vt}{\vt} \in L^1(\Omega)$,
$\frac{1}{\vt} \in L^1(\partial \Omega)$

\item
$\vec{Y}\in W^{1,2}(\Omega)$, $Y_k \geq 0$ a.e. in $\Omega$, $\sumkN Y_k = 1$ a.e. in $\Omega$, $\vf_k\cdot \vc{n}|_{\partial \Omega}=0$

\end{itemize}
satisfy equations (\ref{weak_cont}--\ref{weak_spe}), 
the following entropy inequality
\begin{multline} \label{entropy_ineq}
\int_{\Omega} \frac{ \tn {S} : \nabla \vu}{\vt}\psi \dx
+\int \kappa\frac{|\nabla \vt|^2}{\vt^2}\psi \dx
-\int_{\Omega}\sumkN \omega_k (c_{pk}-c_{vk} \log \vt + \log Y_k)\psi\dx\\
+\int_{\Omega} \psi \sum_{k,l=1}^n D_{kl}\nabla Y_k \cdot \nabla Y_l \dx
+\intpO{\frac{L}{\vt}\vt_0\psi} \leq
\int \frac{\kappa \nabla \vt \cdot \nabla \psi}{\vt} \dx
-\int_{\Omega} \vr s \vu \cdot \nabla \psi \dx\\
-\int_{\Omega} \log \vt \Big(\sumkN \vc{F}_k c_{vk}\Big) \cdot \nabla \psi \dx 
+\int_{\Omega} \Big(\sumkN \vc{F}_k \log Y_k\Big) \cdot\nabla \psi \dx 
+ \intpO{L\psi}
\end{multline}
for all non-negative $\psi \in C^1(\overline{\Omega})$
and the global total energy balance (i.e. (\ref{weak_ene}) with $\psi \equiv 1$)
\begin{equation} \label{glob_ene}
f\intpO{|\vu|^2}+ \intpO{L(\vt-\vt_0)} = \intO{\vr \vc{f}\cdot \vu}.
\end{equation} 
\end{definition}
Formally, the entropy inequality (\ref{entropy_ineq}) is nothing but a weak formulation of the entropy inequality (\ref{24}). However, some modifications are required. First of all, we are not able to keep equality, but due to the technique used to prove existence of such solutions we face the problem that in several terms we are not able to pass to the limit directly and we have to apply the weak lower semicontinuity here.    Note further that  (\ref{entropy_ineq}) does not contain all terms from (\ref{24}), some of them are missing. These terms are formally equal to zero due to assumptions that $\omega_k$ and $\vf_k$ sum up to zero. We removed them from the formulation of the entropy inequality due to the fact that we cannot exclude the situation that $\vr=0$ in some large portions of $\Omega$ (with positive Lebesgue measure), thus $\log \vr$ is not well defined there.  However, the variational entropy solution still has the property that any sufficiently smooth variational entropy solution in the sense above is a classical solution to our problem, provided the density is strictly positive in $\Omega$. Replacing the form of the source terms in the species balance equations by $\vr \omega_k$ we even do not face this problem.  

We are now in position to formulate our main result.

\begin{theorem} \label{t1}
Let $\gamma > 1$, $M>0$, $m > \max\{\frac 23, \frac{2}{3(\gamma-1)}\}$, $a < \frac{3m}{2}$, $\vt_0 \in L^1(\partial \Omega)$, $\vt_0 \geq K_0 >0$ a.e. on $\partial \Omega$.  Let $\Omega \in C^2$ be not axially symmetric. Then there exists at least one  variational entropy solution to Problem P in the sense of Definition \ref{d2}. Moreover, $(\vr,\vu)$ is the renormalized solution to the continuity equation. 

In addition, if $m > 1$, $\gamma > \frac 54$, $a< \frac{3m-2}{2}$, then the solution is a weak solution in the sense of Definition \ref{d1}.

If $\Omega$ is axially symmetric, let $f>0$. Then there exists at least one  variational entropy solution 
to Problem P. In addition, if  $\gamma > \frac 54$, $m > 1$, $m> \frac{16\gamma}{15\gamma -16}$ (if $\gamma \in (\frac 54,\frac 43])$ or $m> \frac{18-6\gamma}{9\gamma -7}$ (if $\gamma \in (\frac 43,\frac 53))$ then the solution is a weak solution. 
\end{theorem}  

\begin{remark} \label{r2}
Recall that the pair $(\vr,\vu)$ is a renormalized solution to the continuity equation provided $\vu \in W^{1,2}(\Omega)$, $\vr \in L^{\frac 65}(\Omega)$ and for any $b \in C^{1}(0,\infty)\cap C([0,\infty))$, $b'(z) = 0$ for $z\geq M$ for some $M>0$
$$
\int_{\Omega} \Big(b(\vr) \Div \psi +  (b(\vr)-b'(\vr)\vr)\Div \vu \psi\Big) \dx =0
$$
for all $\psi \in C^1(\Ov{\Omega})$.
\end{remark}

\section{Proof of the existence results}
    
As explained above, it is enough to prove Theorem \ref{t1} for the case of the generally nondiagonal diffusion matrix. Indeed, for the Fick law, the proof could be simplified due to the special structure of the diffusion flux, but we prefer not to deal with the modification and indicate only one place which is slightly different. The result is exactly the same as in Theorem \ref{t1}.
\begin{proof} { (of Theorem \ref{t1}).} 
First, we define for positive parameters $\delta > \varepsilon > \lambda > \eta >0$ the following approximations of different quantities appearing in the formulation of Problem I. We start with
\begin{equation} \label{101}
{\vc J}_k =-\sumlN Y_k Y_l\widehat{ D}_{kl}(\vt,\vec{Y})\Grad Y_l/Y_l
-\big(\ep(\vr+1) Y_k+\lambda\big)\Grad Y_k/Y_k,\\
\end{equation}
with
\begin{equation} \label{102} 
\widehat {D}_{kl}(\vt,\vec{Y}) = \frac{1}{(\sigma_Y+\ep)^r}  D_{kl}(\vt,\vec{Y})
\end{equation}
for suitably chosen $r\geq 0$,
where $\sigma_Y=\sumkN Y_k$. The reason for this notation is that, unless we let $\lambda \to 0^+$, it is not clear whether $\sigma_Y=1$. We will only know that $Y_k \geq 0$.
For the case of the Fick law this regularization can be simplified. However, in order to keep the unified approach, 
we only slightly modify this step. Instead of (\ref{101}) we set  
\begin{equation} \label{101a}
{\vc J}_k =-\sumlN \widehat{ D}(\vt,\vec{Y})\Grad Y_k
-\big(\ep(\vr+1) Y_k+\lambda\big)\Grad Y_k/Y_k,\\
\end{equation}
where $\widehat D$ is defined similarly as $\widehat{\tn D}$ in (\ref{102}).

Furthermore, we introduce a regularization of the stress tensor 
\begin{equation}\label{103}
\tn{S}_{\eta}={\frac{\mu_{\eta}(\vt)}{1+\eta\vt}}\left[\Grad \vu+(\Grad\vu)^T-\frac{2}{3}\Div \vu \, \tn{I}\right]+{\frac{\nu_\eta(\vt)}{1+\eta\vt}}\lr{\Div \vu}\tn{I},
\end{equation}
where $\mu_{\eta},\nu_{\eta}$ are standard mollifications of the viscosity functions.
Next,
\begin{equation} \label{104}
\kappa_{\delta,\eta} = \kappa^\eta + \delta \vt^B + \delta \vt^{-1}
\end{equation} 
is a regularization of heat conductivity coefficient, where the exponent $B>0$ sufficiently large will
be determined later, and $\kappa^\eta$ is the standard mollification of the heat conductivity.    

We take the following approximation of the specific entropy
\begin{equation} \label{105}
s_k^\lambda=c_{vk}\log \vt - \log Y_k - \log(\vr+\sqrt{\lambda}), 
\end{equation}
and, similarly
\begin{equation} \label{106}
g_k^\lambda = c_{pk} \vt - \vt s_k^\lambda, \quad s^\lambda=\sumkN Y_k s_k^\lambda.
\end{equation}
In what follows, we present only the main steps of the existence proof, pointing always out the specific paper, where more details can be found.

\medskip

\noindent  {\it Step I: Formulation of the approximate problem.}
We consider additionally one more parameter, $N\in \N$ denoting the dimension for the Galerkin approximation of the velocity. 
Let $\{\vc{w}_n\}_{n=1}^{\infty}$ be an orthogonal basis of $W^{1,2}(\Omega)$ such that $\vc{w}\cdot \vc{n}=0$ on $\partial \Omega$ such that 
$\vc{w}_i \in W^{2,q}(\Omega)$ for $q<\infty$ (we can take for example eigenfunctions of the Lam\'e
system with slip boundary conditions).  
We look for
$(\vr_{N,\eta,\lambda,\ep,\delta},\vu_{N,\eta,\lambda,\ep,\delta},\vec{Y}_{N,\eta,\lambda,\ep,\delta},\vt_{N,\eta,\lambda,\ep,\delta})$ 
(from now on we skip the indices) such that\\
$\bullet$ the approximate continuity equation
\begin{equation} \label{108}
\begin{aligned}
\ep\vr+\Div (\vr \vu) &= \ep\Delta\vr+\ep \Ov{\vr},\\
\Grad\vr\cdot\vc{n}\bo&=0,
\end{aligned}
\end{equation}
where $\bar \vr = \frac{M}{|\Omega|}$,
is satisfied pointwisely\\
$\bullet$ the Galerkin approximation for the momentum equation (note that the convective term reduces to the standard form provided $\Div(\vr\vu)=0$, even in the weak sense)
\begin{multline}\label{109}
\intOB{\frac{1}{2}\vr\vu\cdot\Grad\vu\cdot\vc{w}-\frac{1}{2}\vr\left(\vu\otimes\vu\right):\Grad\vc{w}+\tn{S}_{\eta}:\Grad\vc{w}}\\+f \intpO{\vu \cdot \vc{w}}-\intO{(\pi+\delta\vr^{\beta}+\delta\vr^{2})\Div\vc{w}}=\intO{\vr \vc{f}\cdot\vc{w}}
\end{multline}
is satisfied for each test function $\vc{w}\in X_{N}$, where $\vu \in X_N$,
$X_N={\rm span}\{\vc{w}_i\}_{i=1}^N$, and $\beta>0$ is large enough\\ 
$\bullet$ the approximate species mass balance equations 
\begin{equation} \label{110}
\begin{array}{c}
\Div \vc{J}_k=\vw_{k}+\ep \Ov{\vr}_k-\ep Y_k\vr-\Div(Y_k\vr\vu )+\ep\Div(Y_k\Grad\vr) -\sqrt{\lambda} \log Y_k , \\
\vc{J}_k \cdot \vc{n}\bo = 0 
\end{array}
\end{equation}
are satisfied pointwisely,
where $\sumkN \bar \vr_k=\bar \vr$, for example we take $\bar \vr_k=\frac{\bar \vr}{n}$
\\ 
$\bullet$ the approximate internal energy balance
\begin{equation} \label{111}
\begin{aligned}
-\Div\left(\kappa_{\delta,\eta}{\frac{\ep+\vt}{\vt}}\Grad \vt\right)
= &-\Div(\vr e\vu)-\pi\Div\vu +{\frac{\delta}{\vt}} 
+\tn{S}_{\eta}:\Grad{\vu} \\
&+\delta\ep(\beta\vr^{\beta-2}+2)|\Grad\vr|^{2}-\Div\left(\vt\sumkN c_{vk} \vc{J}_k\right)
\end{aligned}
\end{equation}
with the boundary condition
\begin{equation}\label{112}
\kappa_{\delta,\eta} \frac{\ep+\vt}{\vt}\Grad\vt\cdot\vc{n}\bo+(L+\delta\vt^{B-1})(\vt-\vt_{0}^\eta)+\ep \log\vt +\lambda \vt^{\frac B2} \log \vt=0
\end{equation}
is satisfied pointwisely, where $\vt_0^\eta$ is a smooth, strictly positive approximation of $\vt_0$
and $\kappa_{\delta,\eta}$ is as above.

Next, we write down the entropy equality for the approximate system. Note that it is not an additional assumption, but a consequence of the approximate relations above and it is possible to deduce its form (see \cite{PiPo} for more details in the case of the Dirichlet boundary conditions) under the regularity assumptions which correspond to the regularity of solutions to the approximate problem stated above.

\medskip

\noindent {\it Step II: Solvability of the approximate system.}
Following \cite{GPZ} and \cite{PiPo} we can prove
\begin{proposition}\label{p1}
Let $\delta$, $\ep$, $\lambda$ and $\eta$ be positive numbers and $N$ be a positive integer. 
Under the assumptions of Theorem \ref{t1}
there exists a solution to system (\ref{108}--\ref{112}) such that 
$\vr\in W^{2,q}(\Omega)$ $\forall q<\infty$, $\vr\geq0$ in $\Omega$, $\intO{\vr}=M$, $\vu\in X_N$, $\vec{Y}\in W^{1,2}(\Omega)$ with $\log Y_k \in W^{2,q}(\Omega)$ $\forall q<\infty$, $Y_k>0$ a.e. in $\Omega$ and $\vt\in W^{2,q}(\Omega)$ $\forall q<\infty$, $\vt\geq C(N)>0$.
Moreover, this solution satisfies the entropy equation 
\begin{multline} \label{113}
\int_{\Omega} \frac{\psi \tn {S}_\eta : \nabla \vu}{\vt}\dx
+\int_\Omega \kappa_{\delta,\eta}\frac{(\varepsilon+\vt)}{\vt}\frac{|\nabla \vt|^2}{\vt^2}\psi \dx \\
-\int_{\Omega}\omega_k (c_{pk}-c_{vk} \log \vt + \log Y_k)\psi\dx
+\int_{\Omega}\frac{\delta \psi}{\vt^2}\dx \\
-\int_{\Omega} \psi \sumkN \widehat{\vc F}_k \cdot \nabla \log Y_k \dx 
+ \int_{\partial \Omega} \frac{\psi}{\vt}(L+\delta \vt^{B-1})\vt_0^\eta \, \dS \\
+\int_{\Omega}\frac{\delta \varepsilon (\beta \vr^{\beta-2}+2)|\nabla \vr|^2\psi}{\vt}\dx 
+\int_{\Omega}\psi\sumkN(\ep(\vr+1)Y_k+\lambda)\Big|\frac{\nabla Y_k}{Y_k}\Big|^2\dx \\
=\int_{\Omega} \frac{\kappa_{\delta,\eta}(\varepsilon+\vt)\nabla \vt \cdot \nabla \psi}{\vt^2} \dx 
-\int_{\Omega} \vr s^\lambda \vu \cdot \nabla \psi \dx \\
- \int_{\Omega} \sumkN (c_{vk}\log \vt -\log Y_k)  \widehat{\vc F}_k \cdot \nabla \psi \dx - \varepsilon \int_{\Omega} \psi \sumkN Y_k c_{pk} (\Delta \vr + \bar\vr - \vr) \dx\\ 
+ \int_{\Omega} \psi \vr \vu \cdot \Big(\big(\sumkN Y_k\big) \nabla \log (\vr+\sqrt{\lambda}) - \nabla \log \vr\Big)\dx \\ 
+ \int_{\partial \Omega} \frac{\psi}{\vt}\big( (L+\delta\vt^{B-1})\vt + \varepsilon \log \vt + \lambda \vt^{B/2}\log \vt\big)\, \dS \\
-\varepsilon \int_{\Omega} \sumkN Y_k \nabla \vr \cdot \nabla \Big(\frac{g_k^\lambda \psi}{\vt}\Big) \dx
-\sqrt{\lambda} \int_{\Omega} \Big(\sumkN g_k^\lambda  \log Y_k\Big) \frac{\psi}{\vt}\dx \\
+\int_{\Omega}\ep(\Delta \vr+\bar \vr-\vr)(\vr^{\gamma-1}+e+\theta)\frac{\psi}{\vt}\dx\\
+\varepsilon \int_{\Omega} \sumkN (\bar \vr_k - Y_k \vr) \frac{g_k^\lambda \psi}{\vt}\dx 
-\int_{\Omega}\sumkN (\ep(\vr+1)Y_k+\lambda)\frac{\nabla Y_k}{Y_k}\cdot\nabla\psi\dx\\
+\int_{\Omega}\sumkN\big(\ep(\vr+1)Y_k+\lambda\big)s_k^\lambda\frac{\nabla Y_k}{Y_k}\cdot \nabla\psi\dx \\
-\int_{\Omega}\psi\sumkN(\ep(\vr+1)Y_k+\lambda)\frac{\nabla Y_k}{Y_k}\cdot \nabla\log(\vr+\sqrt{\lambda})\dx
\end{multline}
and the following 
estimate
\begin{multline} \label{114}
\sqrt{\lambda}\sumkN \Big(\|Y_k\|_{1,2}+\Big\|\frac{\nabla Y_k}{Y_k}\Big\|_2+\lambda^{-1/4}\|\log Y_k\|_2\Big) + \sumkN \Big\|\frac{|\nabla Y_k|^2}{Y_k} \Big\|_1 
+\|\nabla \vt^{B/2}\|_2 \\ + \Big\|\frac{\nabla \vt}{\vt^{2}}\Big\|_2 
+ \Big\|\frac{\nabla \vr}{\sqrt{\vr+\sqrt{\lambda}}}\Big\|_2
+\|\vt^{-2}\|_1+\|\vt\|_{B,\partial \Omega}+\Big\|\frac{\log \vt}{\vt}\Big\|_{1,\partial \Omega} 
+\|\nabla^2 \vr\|_2 \\ + \|\vu\|_{1,2} + \|\nabla \vr\|_6 \leq C,  
\end{multline}
where $C$ is independent of $N$, $\eta$ and $\lambda$. 
\end{proposition}

Note that if $\Omega$ is axially symmetric (and $f$ is thus positive), the estimate of the $\|\vu\|_2$ (or, more precisely, of $\|\vu\|_{2,\partial \Omega}$) must be deduced from the momentum equation.

\noindent We now let subsequently $N\to +\infty$, $\eta \to 0^+$, $\lambda \to 0^+$, $\varepsilon \to 0^+$, and $\delta \to 0^+$.

\medskip
  
\noindent {\it Step III: Limit passage $N \to +\infty$.}
Using the bounds from Proposition \ref{p1}, weak lower semicontinuity of several terms in the entropy inequality and the fact that for both Galerkin approximation and the limit momentum balance we can use the corresponding velocity as test function (i.e., we have energy equality in both cases), we may let $N\to +\infty$ in the system (\ref{108})--(\ref{113}) above. Note that instead of entropy equality we get entropy inequality.

\medskip

\noindent {\it Step IV: Limit passage $\eta \to 0^+$.}
As we cannot ensure strong convergence of the quadratic term on the rhs of the approximate internal energy balance,
before starting with the limit passage we must replace it by the approximate total energy balance, i.e. we add the kinetic energy balance to the limit version of (\ref{111}). We get
\begin{multline} \label{115}
-\intO{\Big[\vr e+ \frac 12 \vr |\vu|^2 + (\pi+ \delta \vr^\beta + \delta \vr^2)\Big]\vu\cdot \Grad \psi}  
\\
-\intO{\Big(\tn{S}_\eta \vu \cdot \Grad \psi + \delta \vt^{-1} \psi\Big)} 
+\intO{\kappa_{\delta,\eta}{\frac{\ep+\vt}{\vt}}\Grad\vt\cdot\Grad \psi} \\
+ \int_{\partial \Omega}\big[(L+\delta \vt^{B-1})(\vt-\vt_0^\eta) +\ep \log \vt + \lambda \vt^{\frac B2} \log\vt\big] \psi \, \dS + f\intpO{|\vu|^2\psi}\\ 
+ \sum_{k=1}^n c_{vk} \intO{\Big[\vt \sum_{l=1}^n Y_k\widehat{D}_{kl}\Grad Y_l \cdot \Grad \psi + \vt (\ep(\vr+1)Y_k +\lambda) \frac{\Grad Y_k}{Y_k}\cdot \Grad \psi\Big]} \\
 = \intO{\vr \vc{f} \cdot \vu \psi} 
+  \frac{\delta}{\beta-1} \intO{(\ep \beta \Ov{\vr} \vr^{\beta-1}\psi + \vr^\beta \vu \cdot \Grad \psi - \ep \beta \vr^\beta \psi)} \\
+\delta \intO{(2\ep  \Ov{\vr} \vr\psi + \vr^2 \vu \cdot \Grad \psi - 2\ep  \vr^2 \psi)}. 
\end{multline} 
Recalling that the bounds  in (\ref{114}) are independent of $\eta$, it is not difficult to see that we may now let $\eta \to 0^+$ and pass to the limit in our system of equations. Now, if $\Omega$ is axially symmetric, we read the estimate of $\|\vu\|_{2,\partial \Omega}$ from the total energy balance with $\psi$ constant. The same holds also for all subsequent limit passages.

\medskip

\noindent {\it Step V: Limit passage $\lambda \to 0^+$.}
Recall that at this moment it is not yet true that $\vr = \sumkN \vr_k$. However, we have at least
(see (6.12) in \cite{GPZ})
\begin{equation} \label{116}
\|\sumkN \nabla Y_k\|_2 +  \|(\sumkN Y_k)-1\|_6 \leq C(\lambda) \sim \sqrt{\lambda}\to 0 \quad \textrm{for} \quad \lambda \to 0.
\end{equation} 
This bound, together with (\ref{114}), implies
\begin{equation} \label{116a}
\sumkN \|\nabla Y_k\|_{\frac{12}{7}} \leq C
\end{equation} 
with $C$ independent of $\lambda$. We may therefore let $\lambda \to 0^+$ and pass to the limit in our problem (for the details see \cite{PiPo}). Due to (\ref{116}) we see that after the limit passage we have
$$
\sumkN Y_k =1, \qquad \text{ i.e. } \sumkN \vr_k = \vr.
$$

\medskip

\noindent {\it Step VI: Limit passage $\varepsilon \to 0^+$.}
The last two limit passages are nowadays standard in the theory of compressible Navier--Stokes--Fourier system. First, to let $\varepsilon \to 0^+$, we need additional estimates of the total density. We may use the technique of Bogovskii type estimates (see e.g. \cite{NoSt} for more details) to get for $\beta \gg 1$
\begin{equation} \label{117}
\|\vr\|_{\frac 53\beta} \leq C(\delta).
\end{equation}
This can be achieved by testing the approximate momentum equation on the level $\varepsilon >0$ by solution to
$$
\begin{aligned}
\Div \vcg{\vp} &= \vr^{\frac 23 \beta} - \frac{1}{|\Omega|}\intO{\vr^{\frac 23 \beta}}, \\
\vcg{\vp} & = \vc{0} \quad \mbox { on } \partial \Omega.
\end{aligned}
$$
Indeed, this estimate does not imply the compactness of the density sequence and further work must be done: we have to combine the {\it effective viscous flux identity} and the {\it renormalized continuity equation}, see \cite{MPZ} or \cite{NoSt} for more details. After letting $\varepsilon \to 0^+$ we have \\
$\bullet$ the continuity equation
\begin{equation} \label{118}
\int_\Omega \vr \vu\cdot \nabla \psi= 0
\end{equation}
for all $\psi \in C^1(\Ov{\Omega})$ \\
$\bullet$ the weak formulation of the approximate momentum equation
\begin{equation} \label{119}
\begin{array}{c}
\displaystyle
\int_\Omega\Big(-\vr\left(\vu\otimes\vu\right):\Grad\vcg{\varphi}-\tn{S}:\Grad\vcg{\varphi}\Big)\dx+ f\intpO{\vu \cdot \vcg{\vp}}\\
\displaystyle 
-\int_\Omega (\pi+\delta\vr^{\beta}+\delta\vr^{2})\Div\vcg{\varphi} \dx
=\int_\Omega\vr \vc{f}\cdot\vcg{\varphi} \dx
\end{array}
\end{equation}
for all $\vcg{\varphi}\in C^1_{\vc{n}}(\Omega)$\\
$\bullet$ the weak formulation of the approximate species balance equations  
\begin{equation}\label{120}
\begin{array}{c}
\displaystyle  \int_\Omega \Big( -Y_{k}\vr\vu\cdot \Grad \psi +\sum_{l=1}^n Y_k D_{kl}\Grad Y_l \cdot \Grad \psi\Big)\dx  =  \intO{\omega_{k} \psi}
\end{array}
\end{equation}
for all $\psi \in C^1(\Ov{\Omega})$ ($k=1,2,\dots,n$)\\
$\bullet$ the weak formulation of the approximate total energy equation
\begin{equation} \label{121}
\begin{array}{c}
\displaystyle -\intO{\Big[\vr e+ \frac 12 \vr |\vu|^2 + (\pi + \delta \vr^\beta + \delta \vr^2)\Big]\vu\cdot \Grad \psi} 
\\
\displaystyle -\intO{\Big(\tn{S} \vu \cdot \Grad \psi + \delta \vt^{-1} \psi\Big)} 
+\intO{\kappa_{\delta}\Grad\vt\cdot\Grad \psi} + f \intpO{|\vu|^2 \psi}\\
\displaystyle  + \int_{\partial \Omega}\big[(L+\delta \vt^{B-1})(\vt-\vt_0)\big] \psi \, \dS 
\displaystyle + \intO{\vt\sum_{k,l=1}^n c_{vk} Y_k  D_{kl}\Grad Y_l \cdot \Grad \psi} \\
\displaystyle  = \intO{\vr \vc{f} \cdot \vu \psi} + \frac{\delta}{\beta-1} \intO{\vr^\beta \vu \cdot \Grad \psi} 
\displaystyle+\delta \intO{\vr^2 \vu \cdot \Grad \psi} 
\end{array}
\end{equation}
for all $\psi \in C^1(\Ov{\Omega})$ \\
$\bullet$ the weak formulation of the entropy inequality
\begin{multline} \label{122}
\int_{\Omega} \frac{\psi \tn{S} : \nabla \vu}{\vt}\dx
+\int \kappad\frac{|\nabla \vt|^2}{\vt^2}\psi \dx
-\int_{\Omega}\sumkN \omega_k (c_{pk}-c_{vk} \log \vt + \log Y_k)\psi \dx\\
+\int_{\Omega}\frac{\delta \psi}{\vt^2}\dx + \int_{\Omega} \psi \sumkN \sumlN D_{kl}\nabla Y_l \nabla Y_k\dx
+ \int_{\partial \Omega} \frac{\psi}{\vt}(L+\delta \vt^{B-1})\vt_0 \, \dS\\ 
\leq \int \frac{\kappad\nabla \vt \cdot \nabla \psi}{\vt} \dx
-\int_{\Omega} \vr s \vu \cdot \nabla \psi \dx \\
-\int_{\Omega} \sumkN (c_{vk} \log \vt-\log Y_k) \vc{F}_k \cdot\nabla \psi \dx 
+ \int_{\partial \Omega} (L+\delta\vt^{B-1}) \psi \, \dS 
\end{multline}
for all $\psi \in C^1(\Ov{\Omega})$, nonnegative.
More details can be found in \cite{PiPo}. 

\medskip

\noindent {\it Step VII: Estimates independent of $\delta$.} We denote the solution corresponding to $\delta >0$ as $(\vrd,\vud,\vtd,\vec{Y}_{\delta})$.  
The entropy inequality and the total energy balance, both with a constant test function, yield the following estimates for $\Omega$ not axially symmetric
\begin{equation} \label{123}
\|\vtd\|_{1,\partial \Omega} + \delta \|\vtd^B\|_{1,\partial \Omega} \leq C\Big(1+ \Big|\intO{\vrd\vud\cdot \vc{f}}\Big| + \delta \|\vtd^{-1}\|_1\Big).
\end{equation} 
\begin{multline} \label{124}
\|\Grad \vcYd\|_2^2 + \|\Grad \vtd^{\frac m2}\|_2^2 + \|\vud\|_{1,2}^2 + \|\vtd^{-1}\|_{1,\partial \Omega} \\  + \delta \big(\|\Grad \vtd^{\frac B2}\|_2^2 + \|\Grad \vtd^{-\frac 12}\|_2^2 + \|\vtd^{-2}\|_1 +\|\vtd^{B-2}\|_{1,\partial \Omega}\big) \leq C(1+ \delta \|\vtd^{B-1}\|_{1,\partial \Omega}).
\end{multline}
If $\Omega$ is axially symmetric (and $f>0$), then we must add to the left-hand side of (\ref{123}) $\|\vud\|_{2,\partial \Omega}^2$ and replace in the left-hand side of (\ref{124}) $\|\vud\|_{1,2}^2$ by $$
\intO{\frac{\tn{S}(\vtd,\widetilde{\tn{D}}(\vud)):\nabla \vud}{\vtd}}.
$$
Recall also that we know $0\leq (Y_k)_\delta\leq 1$, $k=1,2,\dots, n$.
It is not difficult to bound the $\delta$-dependent terms on the right-hand sides to get for $\Omega$ not axially symmetric
\begin{equation}\label{125}
\begin{aligned}
&\|\Grad \vcYd\|_2+\|\vcYd\|_{\infty} + \|\Grad \vtd^{\frac m2}\|_2 + \|\vud\|_{1,2} + \|\vtd^{-1}\|_{1,\partial \Omega} \\
&+ \delta (\|\Grad \vtd^{\frac B2}\|_2^2 + \|\Grad \vtd^{-\frac 12}\|_2^2 + \|\vtd^{-2}\|_1 +\|\vtd^{B-2}\|_{1,\partial \Omega}) \leq C
\end{aligned}
\end{equation}
and
\begin{equation}\label {vt_3m_0}
\|\vtd\|_{3m} \leq C\Big(1+ \Big|\int_{\Omega}\vrd\vud\cdot\vc{f}\dx\Big| \Big),
\end{equation}
while for $\Omega$ axially symmetric we remove from the left-hand side of (\ref{125}) $\|\vud\|_{1,2}$ and add to the left-hand side of (\ref{vt_3m_0}) $\|\vud\|_{1,2}^2$.
The main issue are now the density estimates. We aim at obtaining 
\begin{equation} \label{136a}
{\rm sup}_{x_0 \in \overline \Omega}\int_{\Omega}\frac{\pi(\vrd,\vtd)+ (1-\alpha)\vrd |\vud|^2}{|x-x_0|^\alpha}\dx \leq C, 
\end{equation}
for some $\alpha>0$ as large as possible. We distinguish three cases. If $x_0$ is sufficiently far from $\partial \Omega$, then it is possible to use as test function in (\ref{119}) 
\begin{equation} \label{129}
\vphi(x) = \frac{x-x_0}{|x-x_0|^\alpha} \tau^2,
\end{equation}
where $\tau$ is a suitable cut-off function. Next, we study the remaining two cases, i.e. $x_0\in \partial \Omega$ and $x_0 \in \Omega$, but close to $\partial \Omega$.  
To simplify the idea, let us assume that we deal with the part of boundary of $\Omega$ which is flat and is described by $x_3 = 0$, i.e. $z(x')=0$, $x' \in {\mathcal O}\subset \R^2$ with the normal vector $\vc{n} = (0,0,-1)$ and  $\vcg{\tau}_1 =(1,0,0)$, $\vcg{\tau}_2=(0,1,0)$ the tangent vectors. The general case can be studied using the standard technique of flattening the boundary, see e.g. \cite{JNP}. Consider  first that $x_0$ lies on the boundary of $\Omega$, i.e. $(x_0)_3 = 0$. Then it is possible to use as the test function in the approximate momentum equation
$$
\vc{w}(x) = \vc{v}(x-x_0),
$$   
where
$$
\vc{v}(x) = 
\frac{1}{|x|^\alpha}(x_1,x_2,x_3) = (x\cdot \vcg{\tau}_1)\vcg{\tau}_1 + (x\cdot \vcg{\tau}_2)\vcg{\tau}_2 + ((0,0,x_3-z(x'))\cdot \vc{n})\vc{n}, \quad x_3\geq 0.
$$
Note that if $(x_0)_3=0$ we get precisely what we need, i.e. estimate (\ref{136a}) (but with $\sup_{x_0 \in \partial \Omega}$ instead of $\sup_{x_0 \in \Ov{\Omega}}$. 
 
However, if $x_0$ is close to the boundary but not on the boundary, i.e. $(x_0)_3>0$, but small, we lose  control of some terms for $0<x_3<(x_0)_3$. In this case, as for the Dirichlet boundary conditions, we must modify the test functions. We first consider
$$
\vc{v}^1(x) = \left\{
\begin{array}{ll}
 \frac{1}{|x-x_0|^\alpha}\big((x-x_0)_1,(x-x_0)_2,(x-x_0)_3\big) , & x_3 \geq \frac{(x_0)_3}{2}, \\[8pt]
 \frac{1}{|x-x_0|^\alpha}\Big((x-x_0)_1,(x-x_0)_2,4(x-x_0)_3 \frac{x_3^2}{|(x-x_0)_3|^2}\Big),  & 0<x_3 <  \frac{(x_0)_3}{2}.
\end{array}
\right.
$$
Nonetheless, using $\vc{v}^1$ as test function we would still miss control of some terms from the convective term, more precisely of those, which contain at least one velocity component $u_3$, however, only close to the boundary, i.e. for $x_3 < (x_0)_3/2$. Hence we further consider
$$
\vc{v}^2(x) = \left\{
\begin{array}{ll}
\displaystyle \frac{(0,0,x_3)}{(x_3+ |x-x_0| |\ln |x-x_0||^{-1})^\alpha} , & |x-x_0| \leq 1/K, \\[8pt]
\displaystyle \frac{(0,0,x_3)}{(x_3+ 1/K |\ln K|^{-1})^\alpha} , & |x-x_0| > 1/K
\end{array}
\right.
$$
for $K$ sufficiently large (but fixed, independently of the distance of $x_0$ from $\partial \Omega$). Note that both functions have zero normal trace, belong to $W^{1,q}(\Omega;\R^3)$ and their norms are bounded uniformly (with respect to the distance of $x_0$ from $\partial \Omega$) provided $1\leq q<\frac 3\alpha$. Thus we finally use as the test function in the approximate momentum balance
\begin{equation} \label{140}
\vcg{\varphi} = \vc{v}^1(x) + K_1 \vc{v}^2(x)
\end{equation}
with $K_1$ suitably chosen (large). Note that the choice of $K$ and $K_1$ is done in such a way that the unpleasant terms from both functions are controlled by those from the other one which provide us a positive information. This is possible due to the fact that the unpleasant terms from $\vc{v}^2$ are multiplied by $|\ln|x-x_0||^{-1} \leq |\ln K|^{-1}\ll 1$. 

We can therefore verify that

\begin{multline}\label{141}
\sup_{x_0 \in \Ov{\Omega}} \intO{ \frac{p(\vrd,\vtd) +\delta(\vrd^\beta + \vrd^2)+ (1-\alpha)\vrd |\vud|^2}{|x-x_0|^\alpha}} \\ 
\leq C (1+ \delta \|\vrd\|_\beta^\beta + \|p(\vrd,\vtd)\|_1 + (1+ \|\vtd\|_{3m})\|\vud\|_{1,2} + \|\vrd|\vud|^2 \|_1),
\end{multline}
provided $0<\alpha < \max\{1, \frac{3m-2}{2m}\}$, $m>\frac 23$,
and, moreover, the test function (\ref{140}) belongs to $W^{1,p}(\Omega;\R^3)$ for $1\leq p < \frac 3\alpha$ with the norm bounded independently of the distance of $x_0$ from $\partial \Omega$. 

We exploit the estimates in the following way. We define now for $1\leq a\leq \gamma$ and $0<b<1$
\begin{equation} \label{II3.11}
\mathcal{B} = \intO{\big(\vrd^a |\vu|^2 + \vrd^b |\vu|^{2b+2}\big)}.
\end{equation}
Then we have
\begin{equation} \label{II3.12}
\|\vrd \vud\|_1 \leq C \mathcal{B}^{\frac{a-b}{2(ab+a-2b)}}.
\end{equation}
and for $1<s< \frac{1}{2-a}$ (if $a<2$), $0<(s-1)\frac{a}{a-1} <b<1$
\begin{equation} \label{II3.13}
\|\vrd |\vud|^2\|_s \leq C \mathcal{B}^{\frac{a-b/s}{ab+a-2b}}.
\end{equation}
Next, we use the Bogovskii-type estimate and get for $1<s<\frac 1{2-a}$ (if $a<2$), $0< (s-1)\frac{a}{a-1} <b<1$, $s\leq \frac{6m}{3m+2}$, $m> \frac 23$
\begin{equation} \label{II3.13a}
\intO{\big(\vrd^{s\gamma} + \vrd^{(s-1)\gamma} p(\vrd, \vtd) + (\vrd |\vud|^2)^s + \delta \vrd^{\beta + (s-1)\gamma}\big)} \leq C(1+ \mathcal{B}^{\frac{sa-b}{ab+a-2b}}).
\end{equation}

We distinguish two cases. First, for $m\geq 2$  the only restriction on $\alpha$ is actually $\alpha <1$. In the other case, if $m \in (\frac 23, 2)$, we have the restriction $\alpha < \frac{3m-2}{2m}$.
Therefore, if $m \geq 2$, we set $a=\gamma$ and using Fatou's lemma and H\"older inequality we show for $b \in ((s-1)\frac{\gamma}{\gamma-1},1)$, $1<s<\frac{2}{2-\gamma}$ and  $s \leq \frac{6m}{3m+2}$ 
\begin{equation} \label{II3.28a}
\begin{array}{c}
\displaystyle 
\sup_{x_0 \in \Ov{\Omega}} \intO{\frac{p(\vrd,\vtd) + (\vrd |\vud|^2)^b}{|x-x_0|}} \\
\displaystyle \leq C \big(1+ \delta \|\vrd\|_\beta^\beta + \|p(\vrd,\vtd)\|_1 + (1+ \|\vtd\|_{3m})\|\vud\|_{1,2} + \|\vrd|\vud|^2 \|_1\big).
\end{array}
\end{equation}

If $m<2$, we keep  $1\leq a<\gamma$ and using H\"older inequality we end up with
\begin{equation} \label{II3.33}
\begin{array}{c}
\displaystyle \sup_{x_0 \in \Ov{\Omega}}\intO{\frac{\vrd^a + (\vrd |\vud|^2)^b}{|x-x_0|}} \\
\displaystyle \leq C \big(1+ \delta \|\vrd\|_\beta^\beta + \|p(\vrd,\vtd)\|_1 + (1+ \|\vtd\|_{3m})\|\vud\|_{1,2} + \|\vrd|\vud|^2 \|_1\big)^{\frac a \gamma} \\
\displaystyle+ C \big(1+ \delta \|\vrd\|_\beta^\beta + \|p(\vrd,\vtd)\|_1 + (1+ \|\vtd\|_{3m})\|\vud\|_{1,2} + \|\vrd|\vud|^2 \|_1\big)^{b}
\end{array}
\end{equation}
which holds for $b \in ((s-1)\frac{\gamma}{\gamma-1},1)$, $1<s<\frac{2}{2-\gamma}$, $\alpha > \max\{\frac{3a-2\gamma}{a}, \frac{3b-2}{b}\}$. 

Let us consider now
\begin{equation} \label{II3.34}
\begin{array}{c}
\displaystyle
-\Delta h = \vrd^a + \vrd^b |\vud|^{2b} - \frac{1}{|\Omega|} \intO{(\vrd^a + \vrd^b |\vud|^{2b})}, \\
\displaystyle
\frac{\partial h}{\partial \vc{n}}|_{\partial \Omega} = 0.
\end{array}
\end{equation}
It is well-known that the unique strong solution admits the following representation
\begin{equation} \label{II3.35}
h(x) = \int_\Omega G(x,y) (\vrd^a + \vrd^b|\vud|^{2b}) \, \mbox{d} y - \frac{1}{|\Omega|} \int_\Omega G(x,y) \, \mbox{d} y \intO{(\vrd^a + \vrd^b |\vud|^{2b})};
\end{equation}
since $G(x,y)\leq C |x-y|^{-1}$, we get due to estimates above
\begin{equation} \label{II3.35a}
\|h\|_\infty \leq C(1+ \mathcal{B}^{\frac{\gamma -b/s}{b\gamma + \gamma -2b}}),
\end{equation}
provided
\begin{equation} \label{II3.35b} 
1<s< \frac{1}{2-\gamma}, \quad 0<(s-1)\frac{\gamma}{\gamma-1} <b<1, \quad s\leq \frac{6m}{3m+2}, \quad m\geq 2,
\end{equation}
and
\begin{equation} \label{II3.35c}
\|h\|_\infty \leq C(1+ \mathcal{B}^{\frac{a -b/s}{ab + a -2b} \frac a \gamma} + \mathcal{B}^{\frac{a -b/s}{ab + a -2b} b}),
\end{equation}
provided
\begin{equation} \label{II3.35d} 
\begin{array}{c}
\displaystyle 1<s< \frac{1}{2-a}, \quad 0<(s-1)\frac{a}{a-1} <b<1, \quad s\leq \frac{6m}{3m+2}, \\[8pt]
\displaystyle \alpha > \frac{3a-2\gamma}{a}, \quad \alpha > \frac{3b-2}{b}, \quad \alpha < \frac{3m-2}{2m}, \quad \frac 23 m\leq 2. 
\end{array}
\end{equation}
Now, from \eqref{II3.11} and \eqref{II3.34}, we have
\begin{equation} \label{II3.36}
\begin{array}{c}
\displaystyle \mathcal{B} = \intO{-\Delta h \vud^2} + \frac{1}{|\Omega|} \intO{\vud^2}\intO{(\vrd^a + \vrd^b |\vud|^{2b})} \\[8pt]
\leq 2\|\nabla \vud\|_2 D^{\frac 12} + C(\varepsilon) \|\vud\|_{1,2}^2 (1+ \mathcal{B}^{\Gamma+\varepsilon})
\end{array}
\end{equation}
for any $\varepsilon >0$, where 
$$ \Gamma =\left\{
\begin{array}{ll}
\displaystyle \frac{\gamma-b}{b\gamma +\gamma -2b} \quad &\mbox{ if } m \geq 2\\
\displaystyle \max\Big\{\frac{a-b}{ab+a-2b} \frac{a}{\gamma}, \frac{a-b}{ab+a-2b} b\Big\} \quad &\mbox{ if } \frac 23 <m<2, 
\end{array}\right.
$$
and
\begin{multline*}
\displaystyle D = \intO{|\nabla h \otimes \vud|^2} -\intO{h \Delta h |\vud|^2} - \intO{h \nabla h \cdot \nabla \vud \cdot \vud} \\
\leq \|h\|_\infty (\mathcal{B} + C(\varepsilon)\|\vud\|_{1,2}^2 \mathcal{B}^{\Gamma +\varepsilon} +\|\nabla \vud\|_2 D^{\frac 12}).
\end{multline*}
Therefore, due to \eqref{II3.35a}, 
\begin{equation} \label{II3.40}
\begin{array}{c}
\mathcal{B} \leq C \big(1+ \mathcal{B}^{\frac{\gamma-b/s}{b\gamma + \gamma-2b}}\big) \quad \mbox{ if } m\geq 2, \\
\mathcal{B} \leq C\big(1+ \mathcal{B}^{\frac {a-b/s}{ab+a-2b}\frac{a}{\gamma}} + \mathcal{B}^{\frac {a-b/s}{ab+a-2b} b}\big) \quad \mbox {if } \frac 23 <m<2.
\end{array}
\end{equation}

Checking carefully all conditions above we end up for $\Omega$ not axially symmetric,
$\gamma >1$ and $m>\frac{2}{4\gamma-3}$, $m>\frac 23$, and for $\Omega$ axially symmetric ($f>0$),  $\gamma >1$ and $m>\frac{6-2\gamma}{3\gamma-1}$, $m>\frac 23$
that there exists $s>1$ such that
\begin{equation} \label{142}
\begin{array}{lcr}
\sup_{\delta>0} \|\vrd\|_{\gamma s} &<& \infty, \\
\sup_{\delta>0} \|\vrd\vud \|_{s} &<& \infty, \\
\sup_{\delta>0} \|\vrd|\vud|^2\|_{s} &<& \infty, \\
\sup_{\delta>0} \|\vud\|_{1,2} &<& \infty, \\
\sup_{\delta>0} \|\vtd\|_{3m} &<& \infty, \\
\sup_{\delta>0} \|\vtd^{m/2}\|_{1,2} &<& \infty, \\
\sup_{\delta>0} \delta \|\vrd^{\beta + (s-1)\gamma}\|_{1} &<& \infty.
\end{array}
\end{equation}
Moreover, we can take $s>\frac 65$ provided 
$\gamma >\frac 54$, $m>\max \{1,\frac{2\gamma +10}{17\gamma -15}\}$ (if $\Omega$ is not axially symmetric),  or $\gamma >\frac 54$, $m> \frac{16\gamma}{15\gamma -16}$ (if $\gamma \in (\frac 54,\frac 43])$ or $m> \frac{18-6\gamma}{9\gamma -7}$ (if $\gamma \in (\frac 43,\frac 53))$ (if $\Omega$ is axially symmetric).

\medskip

\noindent{\it Step VIII: Compactness of the density sequence.}
Here, the standard tools from the theory of compressible Navier--Stokes equations can be applied. We have namely the {\it effective viscous flux identity}
\begin{equation} \label{134}
\begin{array}{c}
\displaystyle
\overline{p(\vr,\vt) T_k(\vr)} - \Big(\frac 43 \mu(\vt) + \xi(\vt)\Big) \overline{T_k(\vr) \Div \,\vu} \\
\displaystyle =  \overline{p(\vr,\vt)} \,\, \overline{T_k(\vr)} - \Big(\frac 43 \mu(\vt) + \xi(\vt)\Big) \overline{T_k(\vr)} \Div \,\vu,
\end{array}
\end{equation}
where
$$
T_k(z) = k T\Big(\frac{z}{k}\Big), \qquad
T(z) = \left\{ \begin{array}{c}
z \mbox{ for } 0\leq z\leq 1, \\
\mbox{ concave on } (0,\infty), \\
2 \mbox{ for } z\geq 3.
\end{array}
\right.
$$  
Next, for the {\it oscillation defect measure}
\begin{equation} \label{135}
\mbox{{\bf osc}}_{\mathbf q} [\vrd\to\vr](Q) = \sup_{k>1} \Big(\limsup_{\delta \to 0^+} \int_Q |T_k(\vrd)-T_k(\vr)|^q \dx\Big)
\end{equation} 
we can show that 
 there exists $q>2$ such that  
\begin{equation} \label{136}
\mbox{{\bf osc}}_{\mathbf{q}} [\vrd\to\vr](\Omega) < \infty,
\end{equation}
provided  $m>\max\{\frac{2}{3(\gamma-1)},\frac 23\}$. This implies the validity of the renormalized continuity equation  even in the case when the density sequence is not bounded in $L^2(\Omega)$. Then it  is not difficult to conclude that
$$
\lim_{k \to \infty} \limsup_{\delta \to 0^+} \intO{|T_k(\vrd)
-T_k(\vr)|^q} = 0
$$
with $q$ as in the estimate of the oscillation defect measure.  Hence
$$
\|\vrd-\vr\|_1 \leq \|\vrd-T_k(\vrd)\|_1 + \|T_k(\vrd) -T_k(\vr)\|_1 + \|T_k(\vr) -\vr\|_1,
$$
which yields strong convergence of the density in $L^1(\Omega)$, and also in $L^p(\Omega)$
for $1 \leq p < s\gamma$. Note that if $s>\frac 65$ we may pass to the limit in the total energy balance while if $s>1$ solely, we may pass to the limit only in the entropy inequality and in the total energy balance with a constant test function. This finishes the proof of Theorem \ref{t1}. 
\end{proof} 

\bibliographystyle{amsalpha}

\end{document}